\documentclass[12pt]{amsart}
\usepackage[centering, hmargin=1in, vmargin=1in]{geometry}
\usepackage{latexsym, amsmath, mathrsfs, bm}
\usepackage[l>=1.1em]{diagrams}
\DeclareMathOperator{\conv}{conv}
\DeclareMathOperator{\eval}{ev}
\DeclareMathOperator{\Ker}{Ker}
\DeclareMathOperator{\Nef}{Nef}
\DeclareMathOperator{\Pic}{Pic}
\DeclareMathOperator{\supp}{Supp}
\DeclareMathOperator{\Sym}{Sym}
\newtheorem{lemma}{Lemma}[section]
\newtheorem{theorem}[lemma]{Theorem}
\newtheorem{corollary}[lemma]{Corollary}
\theoremstyle{definition}
\newtheorem*{conventions}{Conventions}
\newtheorem*{acknowledgements}{Acknowledgements}
\newtheorem{example}[lemma]{Example}
\begin{document}

\title{Syzygies, multigraded regularity and toric varieties}

\author[M.~Hering]{Milena Hering}
\email{mhering@umich.edu}
\address{Department of Mathematics\\University of Michigan\\%
  Ann Arbor\\Michigan 48109\\USA}

\author[H.~Schenck]{Hal Schenck}
\email{schenck@math.tamu.edu} 
\address{Mathematics Department\\Texas A\&M University\\%
  College Station\\Texas 77843\\USA}

\author[G.G.~Smith]{Gregory G. Smith} 
\email{ggsmith@mast.queensu.ca}
\address{Department of Mathematics and Statistics\\%
  Queen's University\\Kingston\\Ontario K7L 3N6\\Canada}


\begin{abstract}
  Using multigraded Castelnuovo-Mumford regularity, we study the
  equations defining a projective embedding of a variety $X$.  Given
  globally generated line bundles $B_{1}, \dotsc, B_{\ell}$ on $X$ and
  $m_{1}, \dotsc, m_{\ell} \in \mathbb{N}$, consider the line bundle
  $L := B_{1}^{m_{1}} \otimes \dotsb \otimes B_{\ell}^{m_{\ell}}$.  We
  give conditions on the $m_{i}$ which guarantee that the ideal of $X$
  in $\mathbb{P}(H^{0}(X,L)^{*})$ is generated by quadrics and the
  first $p$ syzygies are linear.  This yields new results on the
  syzygies of toric varieties and the normality of polytopes.
\end{abstract}

\maketitle

\section{Introduction} 
\label{sec:intro}

Understanding the equations that cut out a projective variety $X$ and
the syzygies among them is a central problem in algebraic geometry.
To give precise statements, consider the morphism $\varphi_{L} \colon
X \rTo \mathbb{P}\bigl( H^{0}(X,L)^{*} \bigr)$ induced by a globally
generated line bundle $L$ on $X$.  Let $S \!=\! \Sym^{\bullet}
H^{0}(X,L)$ be the homogeneous coordinate ring of $\mathbb{P}\bigl(
H^{0}(X,L)^{*} \bigr)$, and let $E_{\bullet}$ be a minimal free graded
resolution of the graded $S$-module $R = \bigoplus_{j \geq 0} H^{0}(X,
L^j)$ associated to $L$.  Following \cite{GreenLazarsfeldI}, we say
that $L$ satisfies \emph{$(N_p)$} for $p \in \mathbb{N}$ provided that
$E_{0} \cong S$ and $E_{i} = \bigoplus S(-i-1)$ for all $1 \leq i \leq
p$.  Thus, $\varphi_{L}(X)$ is projectively normal if and only if $L$
satisfies $(N_{0})$ and $\varphi_{L}(X)$ is normal.  If $L$ satisfies
$(N_{1})$, then the homogeneous ideal of $\varphi_{L}(X)$ is generated
by quadrics and $(N_{2})$ implies that the relations among the
generators are linear.  In \cite{Mumford}, properties $(N_0)$ and
$(N_1)$ are called normal generation and normal presentation
respectively.

Although \cite{GreenII} shows that any sufficiently ample line bundle
on an arbitrary variety satisfies $(N_p)$, it is normally difficult to
determine which multiple of a given ample line bundle suffices.  When
$X$ is a smooth curve of genus $g$, \cite{GreenI} proves that a line
bundle $L$ of degree at least $2g+1+p$ satisfies $(N_p)$.  This is
recovered from an analgous statement for finite sets in
\cite{GreenLazarsfeld}.  When $X$ is a smooth variety of dimension $n$
and $L$ is very ample, \cite{EinLazarsfeld} show that the adjoint line
bundles of the form $K_X + (n+1+p)L$ satisfy $(N_p)$.  Explicit
criteria for $(N_p)$ are also given in \cite{GallegoPurnaprajnaII,
  GallegoPurnaprajna2} for surfaces and in \cite{Pareschi,
  PareschiPopa04} for abelian varieties; we refer to \S1.8.D in
\cite{PAG} for a survey.  The primary goal of this paper is to produce
similar conditions for toric varieties.

To achieve this, we use multigraded Castelnuovo-Mumford regularity.
Fix a list $B_{1}, \dotsc, B_{\ell}$ of globally generated line
bundles on $X$.  For $\bm{u} = (u_1, \dotsc, u_{\ell}) \in
\mathbb{Z}^{\ell}$, set $B^{\bm{u}} := B_{1}^{u_{1}} \otimes \dotsb
\otimes B_{\ell}^{u_{\ell}}$ and let $\bm{\mathcal{B}}$ be the
semigroup $\{ B^{\bm{u}} : \bm{u} \in \mathbb{N}^{\ell} \} \subset
\Pic(X)$.  We say that a line bundle $L$ is
$\mathscr{O}_{X}$\nobreakdash-regular (with respect to $B_{1}, \dotsc,
B_{\ell}$) if $H^{i}(X, L \otimes B^{-\bm{u}}) = 0$ for all $i > 0$
and all $\bm{u} \in \mathbb{N}^{\ell}$ with $|\bm{u}| := u_{1} +
\dotsb + u_{\ell} = i$.  Our main technical result is the following.

\begin{theorem} 
  \label{thm:main}
  Let $\bm{w}_{1}, \bm{w}_{2}, \bm{w}_{3}, \dotsc$ be a sequence in
  $\mathbb{N}^{\ell}$ such that $B^{\bm{w}_{i}} \otimes B_{j}^{-1} \in
  \bm{\mathcal{B}}$ for $1 \leq j \leq \ell$ and set $\bm{m}_{i} :=
  \bm{w}_{1} + \dotsb + \bm{w}_{i}$ for $i \geq 1$.  If
  $B^{\bm{m}_{1}}$ is $\mathscr{O}_{X}$-regular then $B^{\bm{m}_{p}}$
  satisfies $(N_{p})$ for $p \geq 1$.
\end{theorem}

The case $\ell = 1$ is Theorem~1.3 in \cite{GallegoPurnaprajnaII}.
Our proof is a multigraded variant of their arguments.  Applying
Theorem~\ref{thm:main} with $\ell = 1$ to line bundles on toric
varieties yields the following.

\begin{corollary} 
  \label{cor:1}
  Let $L$ be an ample line bundle on an $n$-dimensional toric variety.
  If we have $d \geq n-1+p$ then the line bundle $L^d$ satisfies
  $(N_{p})$.
\end{corollary}

The case $p = 0$, an ingredient in our proof, was established in
\cite{EwaldWessels}; other proofs appear in \cite{LiuTrotterZiegler,
  BrunsGubeladzeTrung, OgataNakagawa}.  On a toric surface,
\cite{Koelman} proves that $L$ satisfies $(N_1)$ if the associated
lattice polygon contains more than three lattice points in its
boundary.  \cite{GallegoPurnaprajna2} gives a criteria for $(N_p)$ on
smooth rational surfaces which, when restricted to toric surfaces,
shows that $L$ satisfies $(N_p)$ if the associated polygon contains at
least $p+3$ lattice points in its boundary.  This result extends to
all toric surfaces and is sharp; see \cite{Schenck}.  Related results
on toric surfaces appear in \cite{Fakhruddin} which studies
multiplication maps and in \cite{Harbourne} which studies $(N_0)$ for
anticanonical rational surfaces.  For an arbitrary toric variety,
\cite{BrunsGubeladzeTrung} shows that $R$ is Koszul when $d \geq n$
and this implies that $L^{d}$ satisfies $(N_{1})$ when $d \geq n$.
Assuming $n \geq 3$, Ogata establishes in \cite{Ogata} that $L^{n-1}$
satisfies $(N_{1})$ and, building on this in \cite{OgataII}, he proves
that $L^{n-2+p}$ satisfies $(N_{p})$ when $n \geq 3$ and $p \geq 1$.

We can strengthen Corollary~\ref{cor:1} by using additional
invariants.  Let
\[
h_{L}(d) := \chi(L^{d}) = \textstyle\sum\limits_{i = 0}^{n} (-1)^{i}
\dim H^{i}(X, L^{d})
\] 
be the Hilbert polynomial of $L$ and let $r(L)$ be the number of
integer roots of $h_{L}$.

\begin{corollary} 
  \label{cor:hilbert}
  Let $L$ be a globally generated line bundle on a toric variety and
  let $r(L)$ be the number of integer roots of its Hilbert polynomial
  $h_{L}$.  If $p \geq 1$ and 
  \[
  d \geq \max\{\deg(h_{L}) - r(L) + p - 1, p\}
  \]
  then the line bundle $L^d$ satisfies $(N_{p})$.
\end{corollary}

If $X = \mathbb{P}^{n}$ and $L = \mathscr{O}_{X}(1)$, then we have
$h_{L}(d) = \binom{d+n}{n}$ and $r(L) = n$.  In particular, we recover
Theorem~2.2 in \cite{GreenII} which states that
$\mathscr{O}_{\mathbb{P}^{n}}(d)$ satisfies $(N_{p})$ for $p \leq d$.
On the other hand, Theorem~2.1 in \cite{OttavianiPaoletti} shows that,
for $n \geq 2$ and $d \geq 3$, $\mathscr{O}_{\mathbb{P}^{n}}(d)$ does
not satisfy $(N_{3d-2})$ and it is conjectured that this is sharp.

Using the dictionary between lattice polytopes and line bundles on
toric varieties, Corollary~\ref{cor:hilbert} yields a normality
criterion for lattice polytopes.  A lattice polytope $P$ is
\emph{normal} if every lattice point in $mP$ for $m \geq 1$ is a sum
of $m$ lattice points in $P$.  Let $r(P)$ be the largest integer such
that $r(P) \, P$ does not contain any lattice points in its interior.

\begin{corollary}
  \label{cor:polytope}
  If $P$ is a lattice polytope of dimension $n$, then
  $\big(n-r(P)\big) \, P$ is normal.
\end{corollary}

Theorem \ref{thm:main} also applies to syzygies of Segre-Veronese
embeddings.

\begin{corollary} 
  \label{cor:prodproj} 
  If $X = \mathbb{P}^{n_{1}} \times \dotsb \times
  \mathbb{P}^{n_{\ell}}$ then the line bundle $\mathscr{O}_{X}(d_{1},
  \dotsc, d_{\ell})$ satisfies $(N_{p})$ for $p \leq \min\{ d_{i} :
  d_i \not= 0 \}$.
\end{corollary}

The Segre embedding $\mathscr{O}_X(1,\dotsc,1)$ satisfies $(N_{p})$ if
and only if $p \leq 3$; see \cite{Lascoux, PragaczWeyman} for $\ell =
2$, and \cite{Rubei, RubeiII} for $\ell > 2$.  An overview of results
and conjectures about the syzygies of Segre-Veronese embeddings
appears in \S3 of \cite{EisenbudGreenHulekPopescu}.

Inspired by \cite{EinLazarsfeld}, we also examine the syzygies of
adjoint bundles.  Recall that a line bundle on a toric variety is
numerically effective (nef) if and only if it is globally generated,
and the dualizing sheaf $K_{X}$ is a line bundle if and only if $X$ is
Gorenstein.

\begin{corollary} 
  \label{cor:canonical}
  Let $X$ be a projective $n$-dimensional Gorenstein toric variety and
  assume that $B_{1}, \dotsc, B_{\ell}$ are the minimal generators of
  $\Nef(X)$.  Suppose $\bm{w}_{1}, \bm{w}_{2}, \dotsc$ is a sequence
  in $\mathbb{N}^{\ell}$ such that $B^{\bm{w}_{i}} \otimes B_{j}^{-1}
  \in \bm{\mathcal{B}}$ for $1 \leq j \leq \ell$ and $\bm{m}_{i} :=
  \bm{w}_{1} + \dotsb + \bm{w}_{i}$ for $i \geq 1$.  If $X \neq
  \mathbb{P}^{n}$ and $p \geq 1$, then the adjoint line bundle $K_{X}
  \otimes B^{\bm{m}_{n+p}}$ satisfies $(N_{p})$.  If $X =
  \mathbb{P}^{n}$ and $p \geq 1$, then $K_{X} \otimes
  B^{\bm{m}_{n+1+p}}$ satisfies $(N_{p})$.
\end{corollary}

\cite{EinLazarsfeld} proves that for a very ample line bundle $L$ and
a nef line bundle $N$ on a smooth $n$\nobreakdash-dimensional
algebraic variety $X \neq \mathbb{P}^n$, $K_X\otimes L^{n+p} \otimes
N$ satisfies $(N_p)$.  Corollary~\ref{cor:canonical} gives a similar
result for ample line bundles on possibly singular Gorenstein toric
varieties.  Specifically, if $L$ is an ample line bundle such that $L
\otimes B_{j}^{-1} \in \bm{\mathcal{B}}$ for $1 \leq j \leq \ell$ and
$N$ is a nef line bundle on $X \neq \mathbb{P}^{n}$ then $K_{X}
\otimes L^{n+p} \otimes N$ satisfies $(N_{p})$.  For an ample line
bundle $L$ on a ruled variety $X$, \S5 in \cite{Butler} demonstrates
that $K_X \otimes L^{n+1+p}$ need not satisfy $(N_p)$ for $p = 0$ and
$1$.  Hence, the conclusions of Corollary~\ref{cor:canonical} are
stronger than one can expect for a general variety.  The proof of
Corollary~\ref{cor:canonical} combines Theorem~\ref{thm:main} with
Fujita's Freeness conjecture for toric varieties, see \cite{Fujino}.

\begin{conventions}
  We work over a field of characteristic zero and $\mathbb{N}$ denotes
  the nonnegative integers.
\end{conventions}

\begin{acknowledgements}
  \begin{flushleft}
  We thank A.~Bayer, W.~Fulton, B.~Harbourne, R.~Lazarsfeld,
  M.~Musta{\c{t}}{\u{a}} and S.~Payne for helpful discussions.  Parts
  of this work were done while the last two authors were visiting the
  Mathematical Sciences Research Institute in Berkeley and the
  Mathematisches Forschungsinstitut in Oberwolfach.  The second author
  was partially supported by NSF Grant DMS 03--11142 and the third
  author was partially supported by NSERC.
  \end{flushleft}
\end{acknowledgements}

\section{Multigraded Castelnuovo-Mumford Regularity} 
\label{sec:reg}

This section reviews multigraded regularity as introduced in
\cite{maclaganSmith}.  Fix a list $B_{1}, \dotsc, B_{\ell}$ of
globally generated line bundles on $X$.  For $\bm{u} := (u_1, \ldots,
u_{\ell}) \in \mathbb{Z}^{\ell}$, set $B^{\bm{u}} := B_{1}^{u_{1}}
\otimes \dotsb \otimes B_{\ell}^{u_{\ell}}$ and let $\bm{\mathcal{B}}$
be the semigroup $\{ B^{\bm{u}} : \bm{u} \in \mathbb{N}^{\ell} \}
\subset \Pic(X)$.  If $\bm{e}_{1}, \dotsc, \bm{e}_{\ell}$ is the
standard basis for $\mathbb{Z}^{\ell}$, then $B^{\bm{e}_{j}} = B_{j}$.

Let $\mathscr{F}$ be a coherent $\mathscr{O}_{X}$-module and let $L$
be a line bundle on $X$.  We say that $\mathscr{F}$ is $L$-regular
(with respect to $B_{1}, \dotsc, B_{\ell}$) provided $H^{i} ( X,
\mathscr{F} \otimes L \otimes B^{- \bm{u}} ) = 0$ for all $i > 0$ and
all $\bm{u} \in \mathbb{N}^{\ell}$ satisfying $|\bm{u}| := u_{1} +
\dotsb + u_{\ell} = i$.  When $X = \mathbb{P}^{n}$, this definition
specializes to Mumford's version of regularity (see
\cite{MumfordCurves}) and as Mumford says, ``this apparently silly
definition reveals itself as follows.''

\begin{theorem} 
  \label{thm:Mumford}
  If the coherent sheaf $\mathscr{F}$ is $L$-regular then for all
  $\bm{u} \in \mathbb{N}^{\ell}$:
  \begin{enumerate}
  \item[(1)] $\mathscr{F}$ is $(L \otimes B^{\bm{u}})$-regular;
  \item[(2)] the map $H^{0} (X, \mathscr{F} \otimes L \otimes
    B^{\bm{u}}) \otimes H^{0}(X, B^{\bm{v}}) \rTo H^{0}(X, \mathscr{F}
    \otimes L \otimes B^{\bm{u}+\bm{v}})$ is surjective for all
    $\bm{v} \in \mathbb{N}^{\ell}$;
  \item[(3)] $\mathscr{F} \otimes L \otimes B^{\bm{u}}$ is globally
    generated provided there exists $\bm{w} \in \mathbb{N}^{\ell}$
    such that $B^{\bm{w}}$ is ample.
  \end{enumerate}
\end{theorem}

\noindent
When $X$ is a toric variety, this follows from results in \S6 of
\cite{maclaganSmith}.  Our approach imitates the proofs of Theorem~2 in
\cite{Mumford} and Proposition~II.1.1. in \cite{Kleiman}.

\begin{proof}
  By replacing $\mathscr{F}$ with $\mathscr{F} \otimes L$, we may
  assume that the coherent sheaf $\mathscr{F}$ is
  $\mathscr{O}_{X}$\nobreakdash-regular.  We proceed by induction on
  $\dim \bigl( \supp (\mathscr{F}) \bigr)$.  The claim is trivial when
  $\dim \bigl( \supp (\mathscr{F}) \bigr) \leq 0$.  As each $B_{j}$ is
  basepoint-free, we may choose a section $s_{j} \in H^{0}(X, B_{j})$
  such that the induced map $\mathscr{F} \otimes B^{- \bm{e}_{j}} \rTo
  \mathscr{F}$ is injective (see page~43 in \cite{Mumford}).  If
  $\mathscr{G}_{j}$ is the cokernel of this map, then we have $0 \rTo
  \mathscr{F} \otimes B^{- \bm{e}_{j}} \rTo \mathscr{F} \rTo
  \mathscr{G}_{j} \rTo 0$ and $\dim \bigl( \supp(\mathscr{G}_{j})
  \bigr) < \dim \bigl( \supp (\mathscr{F}) \bigr)$.  From this short
  exact sequence, we obtain the long exact sequence
  \[
  H^{i}(X, \mathscr{F} \otimes B^{- \bm{u} - \bm{e}_{j}}) \rTo
  H^{i}(X, \mathscr{F} \otimes B^{- \bm{u}}) \rTo H^{i}(X,
  \mathscr{G}_{j} \otimes B^{- \bm{u}}) \rTo H^{i+1}(X, \mathscr{F}
  \otimes B^{- \bm{u}- \bm{e}_{j}}) \rTo \, .
  \]
  By taking $|\bm{u}| = i$, we deduce that $\mathscr{G}_{j}$ is
  $\mathscr{O}_{X}$-regular. The induction hypothesis implies that
  $\mathscr{G}_{j}$ is $(B_{j})$-regular.  Setting $\bm{u} =
  - \bm{e}_{j} + \bm{u}'$ with $|\bm{u}'| = i$, we see that
  $\mathscr{F}$ is $(B_{j})$-regular and (1) follows.

  For (2), consider the commutative diagram:
  \[
  \begin{diagram}[vtrianglewidth=1em]
    &&&& H^{0}(X, \mathscr{F}) \otimes H^{0}(X,B_{j}) & \rTo &
    H^{0}(X, \mathscr{G}_{j}) \otimes H^{0}(X, B_{j}) \\
    &&& \ruTo(2,2) & \dTo && \dTo \\
    0 &\rTo & H^{0}(X, \mathscr{F}) & \rTo & H^{0}(X, \mathscr{F}
    \otimes B_{j}) & \rTo & H^{0}(X, \mathscr{G}_{j}\otimes B_{j}) \,
    .
  \end{diagram}
  \]
  Since $\mathscr{F}$ is $\mathscr{O}_{X}$-regular, the map in the top
  row is surjective.  The induction hypothesis guarantees that the map
  in the right column is surjective.  Thus, the Snake Lemma
  (e.g. Prop.~1.2 in \cite{Bourbaki}) implies that the map in the
  middle column is also surjective.  Therefore, (2) follows from the
  associativity of the tensor product and (1).

  Lastly, consider the commutative diagram:
  \[
  \begin{diagram}[small]
    H^{0}(X, \mathscr{F} \otimes B^{\bm{u}}) \otimes H^{0}(X,
    B^{\bm{v}}) \otimes \mathscr{O}_{X} & \rTo & H^{0}(X, \mathscr{F}
    \otimes B^{\bm{u+v}}) \otimes \mathscr{O}_{X} \\
    \dTo && \dTo^{\beta_{\bm{u} + \bm{v}}}\\
    H^{0}(X, \mathscr{F} \otimes B^{\bm{u}}) \otimes B^{\bm{v}} &
    \rTo^{\beta_{\bm{u}}\otimes \text{id}} & \mathscr{F} \otimes
    B^{\bm{u} + \bm{v}}
  \end{diagram}
  \]
  Applying (2), we see that the map in the top row is surjective.  By
  assumption, there is $\bm{w} \in \mathbb{N}^{\ell}$ such that
  $B^{\bm{w}}$ is ample.  If $\bm{v} := k \bm{w}$, then Serre's
  Vanishing Theorem (e.g. Theorem~1.2.6 in \cite{PAG}) implies that
  $\beta_{\bm{u} + \bm{v}}$ is surjective for $k \gg 0$.  Hence,
  $\beta_{\bm{u}}$ is also surjective which proves (3).
\end{proof}

We end this section with an elementary observation.

\begin{lemma} 
  \label{lem:ses} 
  Let $0 \rTo \mathscr{F}' \rTo \mathscr{F} \rTo \mathscr{F}'' \rTo 0$
  be a short exact sequence of coherent
  $\mathscr{O}_{X}$\nobreakdash-modules.  If $\mathscr{F}$ is
  $L$-regular, $H^{0}(X, \mathscr{F} \otimes L \otimes B^{-
    \bm{e}_{j}}) \rTo H^{0}(X, \mathscr{F}'' \otimes L \otimes B^{-
    \bm{e}_{j}})$ is surjective for all $1 \leq j \leq \ell$, and
  $\mathscr{F}''$ is $(L \otimes B^{-\bm{e}_{j}})$-regular for all $1
  \leq j \leq \ell$, then $\mathscr{F}'$ is also $L$-regular.
\end{lemma}

\begin{proof}[Sketch of Proof]
  Tensor the exact sequence $0 \rTo \mathscr{F}' \rTo \mathscr{F} \rTo
  \mathscr{F}'' \rTo 0$ with $L \otimes B^{\bm{u}}$ and analyze the
  associated long exact sequence.  The argument is similar to the
  proof of Theorem~\ref{thm:Mumford}.1.
\end{proof}

\section{Proof of Main Theorem} 
\label{sec:main}

The proof of Theorem~\ref{thm:main} combines multigraded
Castelnuovo-Mumford regularity with a cohomological criterion for
$(N_{p})$.  Given a globally generated line bundle $L$ on $X$, there
is a natural surjective map $\eval_{L} \colon H^{0}(X, L) \otimes
\mathscr{O}_{X} \rTo L$ and we set $M_{L} := \Ker(\eval_{L})$.  Hence,
$M_{L}$ is a vector bundle on $X$ which sits in the short exact
sequence
\begin{equation} \tag{$\dagger$}
  0\rTo M_{L} \rTo H^0(X,L)\otimes \mathscr{O}_X \rTo L \rTo 0 \, .
\end{equation}
It is well-known that $M_{L}$ governs the syzygies of $\varphi_{L}(X)$
in $\mathbb{P} \bigl( H^{0}(X, L)^{*} \bigr)$.  Specifically, $L$
satisfies $(N_{p})$ if and only if $H^1(X, \bigwedge^{q} M_{L} \otimes
L^j) = 0$ for $q \leq p+1$ and $j \geq 1$; see Lemma~1.10 in
\cite{GreenLazarsfeld} or Proposition~1.3.3 in \cite{Lazarsfeld}.  In
characteristic zero, $\bigwedge^{k} M_{L}$ is a direct summand of
$M_{L}^{\otimes k}$, so it suffices to show that $H^{1}(X,
M_{L}^{\otimes q} \otimes L^{j}) = 0$ for $q \leq p+1$ and $j \geq 1$
in our situation.

\begin{proof}[Proof of Theorem~\ref{thm:main}]
  Set $L := B^{\bm{m}_{p}}$ and let $M_{L}$ be the vector bundle in
  $(\dagger)$.  We first prove, by induction on $q$, that
  $M_{L}^{\otimes q}$ is $(B^{\bm{m}_{q}})$-regular for all $q \geq
  1$.  Since $B^{\bm{m}_{1}}$ is $\mathscr{O}_{X}$-regular,
  Theorem~\ref{thm:Mumford}.2 implies that $H^{0}(X,B^{\bm{m}_{1} +
    \bm{u}}) \otimes H^{0}(X,B^{\bm{v}}) \rTo H^{0}(X, B^{\bm{m}_{1} +
    \bm{u} + \bm{v}})$ is surjective for all $\bm{u}, \bm{v} \in
  \mathbb{N}^{\ell}$.  In particular, the maps $H^{0}(X,L) \otimes
  H^{0}(X, B^{\bm{m}_{1} - \bm{e}_{j}}) \rTo H^{0}(X,L \otimes
  B^{\bm{m}_{1} - \bm{e}_{j}})$ for $1 \leq j \leq \ell$ are
  surjective because $B^{\bm{m}_{1}} \in \bigcap_{j = 1}^{\ell} (
  B_{j} \otimes \bm{\mathcal{B}})$.  Combining
  Theorem~\ref{thm:Mumford}.1 and Lemma~\ref{lem:ses}, we see that
  $M_L$ is $(B^{\bm{m}_{1}})$\nobreakdash-regular.  For $q > 1$,
  tensor the sequence $(\dagger)$ with $M_{L}^{\otimes (q-1)}$ to
  obtain the exact sequence $0 \rTo M_{L}^{\otimes q} \rTo H^{0}(X,L)
  \otimes M_{L}^{\otimes (q-1)} \rTo M_{L}^{\otimes (q-1)} \otimes L
  \rTo 0$.  The induction hypothesis states that $M_{L}^{\otimes
    (q-1)}$ is $(B^{\bm{m}_{q-1}})$-regular.  Since $B^{\bm{w}_{q}}
  \otimes B_{j}^{-1} \in \bm{\mathcal{B}}$ for all $1 \leq j \leq
  \ell$, Theorem~\ref{thm:Mumford}.2 shows that 
  \[
  H^{0}(X, M_{L}^{\otimes (q-1)} \otimes B^{\bm{m}_{q} - \bm{e}_{j}})
  \otimes H^{0}(X, L) \rTo H^{0}(X, M_{L}^{\otimes (q-1)} \otimes L
  \otimes B^{\bm{m}_{q} - \bm{e}_{j}})
  \] 
  is surjective for $1 \leq j \leq \ell$.  Again by
  Theorem~\ref{thm:Mumford}.1 and Lemma~\ref{lem:ses}, $M_{L}^{q}$ is
  $(B^{\bm{m}_{q}})$-regular.

  As observed above, it suffices to prove that $H^{1}(X,
  M_{L}^{\otimes q} \otimes L^{j}) = 0$ for $q \leq p + 1$ and $j \geq
  1$.  Since $M_{L}^{\otimes q}$ is $(B^{\bm{m}_{q}})$-regular,
  Theorem~\ref{thm:Mumford}.1 implies that $M_{L}^{\otimes q}$ is
  $(B^{\bm{m}_{p}})$-regular for $1 \leq q \leq p$; as
  $\mathscr{O}_{X}$ is $(B^{\bm{m}_{1}})$-regular,
  Theorem~\ref{thm:Mumford}.1 also implies that $\mathscr{O}_{X}$ is
  $(B^{\bm{m}_{p}})$-regular.  It follows that $H^{1}(X,
  M_{L}^{\otimes q} \otimes L^{j}) = 0$ for $q \leq p$ and $j \geq 1$.
  Moreover, Theorem~\ref{thm:Mumford}.2 shows that
  $H^{0}(X,L) \otimes H^{0}(X, M_{L}^{\otimes p} \otimes L^{j}) \rTo
  H^{0}(X, M_{L}^{\otimes p} \otimes L^{j+1})$
  is surjective and the exact sequence $0 \rTo M_{L}^{\otimes(p+1)}
  \otimes L^{j} \rTo H^{0}(X,L) \otimes M_{L}^{\otimes p} \otimes L^{j}
  \rTo M_{L}^{\otimes p} \otimes L^{j+1} \rTo 0$ implies that $H^{1}(X,
  M_{L}^{\otimes(p+1)} \otimes L^{j}) = 0$ for $j \geq 1$.
\end{proof}

\section{Applications to Toric Varieties}
\label{sec:toric}

In this section, we apply the main theorem to line bundles on an
$n$-dimensional projective toric variety $X$.  Consider a globally
generated line bundle $L$ on $X$ and its associated lattice polytope
$P_{L}$.  Let $r(L)$ be the number of integer roots of the Hilbert
polynomial $h_{L}(d) = \chi(L^d)$.  Since the higher cohomology of
$L^{d}$ vanishes and the lattice points in the polytope $d P_{L} =
P_{L^{d}}$ form a basis for $H^{0}(X,L^{d})$, it follows that
$h_{L}(d)$ equals the Ehrhart polynomial of $P_{L}$.  In other words,
$h_{L}(d)$ is the number of lattice points in $dP$.  If $r(P_{L})$ is
the largest integer such that $r(P_{L}) P_{L}$ does not contain any
interior lattice points, then Ehrhart reciprocity
(e.g. Corollary~4.6.28 in \cite{Stanley}) implies that the integer
roots of $h_{L}(d)$ are $\{-1, \dotsc , -r(P_{L}) \}$ and $r(P_{L}) =
r(L)$.

\begin{lemma} 
  \label{lem:regularity}
  If $L$ is a globally generated line bundle on a toric variety $X$
  and $r(L)$ is the number of integer roots of its Hilbert polynomial
  $h_{L}$, then $L^{\deg(h_L)-r(L)}$ is $\mathscr{O}_{X}$-regular with
  respect to $L$.
\end{lemma}

\begin{proof}
  We must prove that $H^{i}(X, L^{\deg(h_L) - r(L)-i}) = 0$ for all $i
  > 0$.  If $\deg(h_L) - r(L) - i \geq 0$, this follows from the
  vanishing of the higher cohomology of globally generated line
  bundles on a complete toric variety; see \S3.5 of \cite{Fulton}.
  When $\deg(h_L) - r(L) - i < 0$, we follow the proof of Theorem~2.5
  in \cite{BatyrevBorisov}.  Let $X'$ be the toric variety
  corresponding to the normal fan to $P_L$.  There is a canonical
  toric map $\psi \colon X \rTo X'$ and an ample line bundle $A$ on
  $X'$ such that $H^{i}(X, L^r) \cong H^{i}(X', A^r)$ for all $r \in
  \mathbb{Z}$.  A toric version of the Kodaira Vanishing Theorem
  establishes that $H^{j}(X, L^{-u}) = 0$ for $u >0$ and $j \neq
  \deg(h_{L}) = \dim P_L = \dim X'$ (combine Serre duality from \S4.4
  of \cite{Fulton} with Theorem~3.4 in \cite{Mustata}). In particular,
  we have $H^{i}(X, L^{\deg(h_L) - r(L) - i}) = 0$ for $i \neq
  \deg(h_L)$.  When $i = \deg(h_L)$, we also have $0 = h_{L}\bigl(
  -r(L) \bigr) = \chi(L^{-r(L)}) = (-1)^{i} \dim H^{i}(X,
  L^{\deg(h_{L}) -r(L) -i})$.
\end{proof}

\begin{proof}[Proof of Corollary~\ref{cor:hilbert}]
  In light of Lemma~\ref{lem:regularity}, the claim is the special
  case of Theorem~\ref{thm:main} with $\ell = 1$, $B_{1} = L$,
  $\bm{w}_{1} = \max\{\deg(h_L) - r(L),1\}$ and $\bm{w}_{i} = 1$ for
  $i > 1$.
\end{proof}

\begin{proof}[Proof of Corollary~\ref{cor:1}]
  The case $p = 0$ is in \cite{EwaldWessels}; the case $p \geq 1$
  follows from Corollary~\ref{cor:hilbert}.
\end{proof}

The following well-known examples illustrate that
Corollary~\ref{cor:hilbert} is sharp in some cases.

\begin{example}
  \label{exa}
  Let $L$ be the ample line bundle on the toric variety $X$
  corresponding to the polytope $\conv\{(1,0), (0,1), (1,1), (2,2) \}
  \subset \mathbb{R}^2$.  The ideal of $X$ in $\mathbb{P}^{3} =
  \mathbb{P}\bigl( H^0(X,L)^{*} \bigr)$ is generated by the cubic
  $x_2^3 - x_0 x_1x_3$ which implies that $L$ does not satisfy
  $(N_1)$.  Calculations in \cite{M2} show that $L^2$ satisfies
  $(N_3)$ but not $(N_4)$.
\end{example}

\begin{example}
  Let $\bm{e}_{1}, \dotsc, \bm{e}_{n}$ be the standard basis of
  $\mathbb{R}^{n}$ and let $L$ be the ample line bundle on the toric
  variety $X$ corresponding to 
  \[
  P = \conv\{\bm{0}, \bm{e}_{1}, \dotsc, \bm{e}_{n-1}, \bm{e}_{1} +
  \dotsb + \bm{e}_{n-1} + (n-1)\bm{e}_{n} \} \subset \mathbb{R}^n \, .
  \]
  Since $X$ is $n$-dimensional and singular, the morphism $X \rTo
  \mathbb{P}^{n} = \mathbb{P}\bigl( H^{0}(X,L)^* \bigr)$ is obviously
  not an embedding.  The natural map $S \rTo R$ is not surjective and
  $(n-2)P$ is not normal because $(1, \dotsc, 1)$ lies in $2(n-2)P$
  but cannot be written as an integral combination of two lattice
  points in $(n-2)P$.  For $n = 3$, calculations in \cite{M2} show
  that $L^2$ satisfies $(N_1)$ but not $(N_2)$.
\end{example}

\begin{proof}[Proof of Corollary~\ref{cor:polytope}]
  Given a lattice polytope $P$, let $X$ be the corresponding toric
  variety and $L$ the associated ample line bundle on $X$.  Since $P$
  is normal if and only if $L$ satisfies $(N_0)$, the result follows
  from Corollary~\ref{cor:hilbert} and the fact that $r(P) = r(L)$.
\end{proof}

\begin{proof}[Proof of Corollary~\ref{cor:prodproj}]
  Let $\pi_{i} \colon X \rTo \mathbb{P}^{n_i}$ be the $i$-th
  projection and set $B_{i} := \pi_{i}^{*}\bigl(
  \mathscr{O}_{\mathbb{P}^{n_{i}}}(1) \bigr)$.  If $I := \bigl\{ i \in
  \{1, \dotsc, \ell \} : d_{i} \neq 0 \bigr\}$, then $\mathscr{O}_{X}
  (d_{1}, \dotsc, d_{\ell}) \cong \bigotimes_{i\in I} B_{i}^{d_{i}}$.
  Let $\bm{\mathcal{B}}$ be the semigroup generated by $\{ B_i : i \in
  I \}$ and let $d := \min\{d_{i}-1 : i \in I\}$.  Proposition~6.10 in
  \cite{maclaganSmith} proves that $\mathscr{O}_{X}$ is
  $\mathscr{O}_{X}$-regular with respect to $B_{1}, \dotsc, B_{\ell}$.
  Thus, Theorem~\ref{thm:Mumford} shows that $\bigotimes_{i\in I}
  B_{i}^{d_{i}-d} $ is $\mathscr{O}_{X}$\nobreakdash-regular with
  respect to $\{B_i: i \in I\}$ and lies in $\bigcap_{i \in I} ( B_{i}
  \otimes \bm{\mathcal{B}})$.  Since we have $\bigotimes_{i \in I} B_i
  \in \bigcap_{j\in I} ( B_{j} \otimes \bm{\mathcal{B}})$,
  Theorem~\ref{thm:main} applies with $\bm{w}_{1} = (d_{1}-d, \dotsc,
  d_{\ell} - d)$ and $\bm{w}_{j} = (1, \dotsc, 1)$ for $j \geq 2$.
\end{proof}

Now assume that $B_1, \dotsc, B_\ell$ are the minimal generators of
$\Nef(X)$.  To apply our techniques to adjoint bundles, we need to
find $\bm{u}$ with $K_{X} \otimes B^{\bm{u}} \in \bm{\mathcal{B}} =
\Nef(X)$.  Inspired by Fujita's conjectures, Corollary~0.2 in
\cite{Fujino} provides the following necessary criterion: \emph{Let
  $X$ be a projective toric variety $($not isomorphic to
  $\mathbb{P}^n)$ such that the canonical divisor $K_X$ is
  $\mathbb{Q}$-Cartier.  If $D$ is a $\mathbb{Q}$-Cartier divisor such
  that $D \cdot C \geq n$ for all torus invariant curves $C$, then
  $K_{X} + D$ is nef.}
 
\begin{proof}[Proof of Corollary~\ref{cor:canonical}]
  If $X = \mathbb{P}^{n}$, then $K_{X} = \mathscr{O}_{X}(-n-1)$;
  Corollary~\ref{cor:hilbert} proves that $K_{X} \otimes
  B^{\bm{m}_{n+1+p}}$ satisfies $(N_{p})$.  Theorem~3.4 in
  \cite{Mustata} shows that $K_{X} \otimes B^{\bm{m}_{n+1}}$ is
  $\mathscr{O}_{X}$-regular with respect to $B_{1}, \dotsc, B_{\ell}$.
  For any torus invariant curve $C$, there is a $B_{i}$ such that
  $B_{i} \cdot C> 0$.  Since $B^{\bm{m}_{n}} = B_{i}^{n} \otimes B'$
  where $B'$ is globally generated, Corollary~0.2 in \cite{Fujino}
  implies that $K_{X} \otimes B^{\bm{m}_{n}} \in \bm{\mathcal{B}}$.
  It follows that $K_{X} \otimes B^{\bm{m}_{n+1}} \in \bigcap_{j =
    1}^{\ell} ( B_{j} \otimes \bm{\mathcal{B}})$ and
  Theorem~\ref{thm:main} proves the claim.
\end{proof}

The singular cubic surface in Example~\ref{exa} also demonstrates that
Corollary~\ref{cor:canonical} can be sharp; see
\cite{GallegoPurnaprajnaII} for more examples of this type.

\begin{example}
  Let $(X,L)$ be the normal cubic surface and ample line bundle
  defined in Example~\ref{exa}.  It follows that $K_X^{-1} = L$ and
  $L$ is the minimal generator of $\Nef(X)$.  Example~\ref{exa} shows
  that $K_{X} \otimes L^2 = L$ does not satisfy $(N_1)$.  Hence,
  Corollary~\ref{cor:canonical} provides the smallest $m \in
  \mathbb{N}$ (namely $m = 3$) such that $K_X \otimes L^{m}$ satisfies
  $(N_1)$.
\end{example}

For toric surfaces, it follows from \cite{Schenck} that all of our
corollaries are not optimal for $p \geq 2$.  Specifically, given an
ample line bundle $L \not\cong \mathscr{O}_{\mathbb{P}^2}(1)$ on a
Gorenstein toric surface $X$, $K_X \otimes L^m$ satisfies $N_{3(m-2)}$
for $m \geq 2$ and $m \neq 4$, and $K_X \otimes L^4$ satisfies $N_5$.

\raggedright
\providecommand{\bysame}{\leavevmode\hbox 
  to3em{\hrulefill}\thinspace}

\end{document}